\documentclass[11pt]{amsart}
\numberwithin{equation}{section}
\newtheorem{theorem}{Theorem}

\numberwithin{theorem}{section}
\newtheorem{remark}[theorem]{Remark}

\def\al{\aligned}
\def\eal{\endaligned}
\def\be{\begin{equation}}
\def\ee{\end{equation}}
\def\lab{\label}

\def\M{{\bf M}}

\def\al{\aligned}

\numberwithin{equation}{section}

\begin{document}

\tracingpages 1
\title[diameter bound]
{\bf On the question of diameter bounds in Ricci flow}
\author{ Qi S. Zhang}
\address{
Department of
Mathematics, University of California, Riverside, CA 92521, USA}
\date{\today;  MSC 	53C44,  	35K40,  	53C20 }

\begin{abstract}
A question about Ricci flow is when the diameters of the manifold under the evolving metrics stay finite and bounded away from $0$.
 Topping \cite{T:1} addresses the question with an upper bound
that depends on the $L^{(n-1)/2}$ bound of the scalar curvature, volume and a local version of Perelman's $\nu$ invariant. Here $n$ is the dimension.
His result is sharp when Perelman's F entropy is positive. In this note, we give a direct
proof that for all compact manifolds,
the diameter bound depends just on
the $L^{(n-1)/2}$ bound of the scalar curvature, volume and the Sobolev constants (or positive Yamabe constant).
This bound seems directly computable in large time for some Ricci flows.
In addition, since the result in its most general form is independent of Ricci
flow, further applications may be possible.

A generally sharp lower bound for the diameters is also given, which depends only
on the initial metric, time and $L^\infty$ bound of the scalar curvature.
These results imply that, in finite time, the Ricci flow can neither turn the diameter to infinity nor zero, unless the scalar curvature blows up.
 \end{abstract}
\maketitle
\section{statement of result}

The Ricci flow introduced by R. Hamilton is a nonlinear parabolic
equation along which the metrics of a Riemannian manifold evolve.
Therefore understanding the evolution of basic geometric quantities
such as volumes, diameters etc has been a basic task in the study of
the Ricci flow. By now it is known that in finite time, the volume
of geodesic balls along the Ricci flow are comparable to that of the
Euclidean balls provided that the scalar curvatures are bounded.
These are the $\kappa$ non-collapsing property by Perelman and the
so called $\kappa$ non-inflating property. Perelman also proved that
the diameter stays bounded for the (normalized) K\"ahler Ricci flow
on Fano manifolds. See the paper by Sesum and Tian \cite{ST:1}.
Ilmanen and Knopf \cite{IK:1} proved a lower bound for the diameters
under a topological condition. In
the paper \cite{T:1}, P. Topping considered the question of diameter upper
bounds along a general Ricci flow $\partial_t g = -2 Ric $ on a compact Riemannian manifold $M$ of dimension $n$. One of his main result (Theorem 2.4) says that  there exists a constant $C=C(g(0), T)$ such that for all $t \in [0, T)$,
if $diam (M, g(t)) \ge C$, then
\[
diam (M, g(t)) \le C \int_M R^{(n-1)/2} dg(t).
\] Here $R$ is the scalar curvature. The constant $C$ depends on volume and
a local version of Perelman's $\nu$ entropy.
If the infimum of the $F$ entropy for $(M, g(0))$ is positive,
then the above bound holds without the lower bound assumption; and
$C$ is independent of time.
See also Theorem 6.75 in  B. Chow etc \cite{C++:1} for
another exposition of this result. Let us recall that Perelman's
$F$ entropy is $F(v)=\int_M (4 |\nabla v|^2 + R v^2) dg$ where $v$
is a $W^{1, 2}$ function on $M$ with unit $L^2$ norm.

The first goal
of this note  is to prove that:
for all compact manifolds,
the diameter bound depends just on
the $L^{(n-1)/2}$ bound of the scalar curvature, volume and the Sobolev constants to be defined below. The Sobolev constant can also be replaced by the
positive Yamabe constant. This result is similar in spirit to the classical Bonnet- Myers  which says that a positive lower bound of the Ricci curvature implies
an upper bound of the diameter.  The difference is that the Ricci lower bound
is replaced the three quantities mentioned above. Note that when $n=3$,
the integral involving the scalar curvature is just the total curvature.
We mention that some relation
between the Yamabe constant and the area of black holes are found in the
paper by Cai and Galloway \cite{CG:1}.

Another goal of the note is to prove a lower bound for the diameters, which depends only
on the initial metric, time and $L^\infty$ bounds of the scalar curvature. Therefore, if the scalar curvature is bounded,
the diameter of the manifold at time $t$ is comparable to that
of the initial manifold.

To prove the upper bound, we will
build on the idea in \cite{T:1} where, Perelman's $W$ entropy, a
maximal type of function for $R_+$ and a covering technique  are
used. The new input is the uniform Sobolev inequality along Ricci
flow and a partition and covering argument involving the
volume of the manifold.
This bound seems directly computable in large time for some Ricci flows.
 This is the case when $M$ is a compact quotient
of the hyperbolic space. For these manifolds, the Sobolev constants are independent of time, since
the Ricci curvature and injectivity radii are bounded from below uniformly.
In addition, since the result in its most general form is independent of Ricci
flow, further applications may be possible.

Now let us define the Sobolev constants mentioned above.
Let $(M, g)$ be a $n$ dimensional compact Riemannian manifold and $R$
be the scalar curvature. It is well known that the following Sobolev inequality holds.  For any $v \in W^{1, 2}(M)$, there exist positive constants
$A$ and $B$, depending on $g$ such that
\be
\lab{sob}
\left(\int_M v^{2n/(n-2)} dg\right)^{(n-2)/n} \le
A \int_M ( 4 | \nabla v|^2 + R v^2 ) dg + B \int_M v^2 dg.
\ee When $(M, g(t))$ is a Ricci flow, the following uniform Sobolev inequality
also holds. See Theorem 6.2.1
in \cite{Z11:1} e.g..
For any $v \in W^{1, 2}(M)$, there exist positive constants
$A$ and $B$, depending only on the Sobolev constant of $g(0)$, $R(\cdot, 0)$ and $t$ such that
\be
\lab{sobRF}
\left(\int_M v^{2n/(n-2)} dg(t)\right)^{(n-2)/n} \le
A \int_M ( 4 | \nabla v|^2 + R v^2 ) dg(t) + B \int_M v^2 dg(t).
\ee Also, if the infimum of the $F$ entropy is positive, then $B=0$ and
$A$ is independent of time.

The Yamabe constant is
\be
\lab{yama}
Y(g)= \inf_{v \in C^\infty(M), v>0}
\frac{\int_M \left(4 (n-1) (n-2)^{-1} |\nabla v|^2 + R v^2 \right) dg}{\left(
\int_M v^{2n/(n-2)} dg \right)^{(n-2)/n}}.
\ee

Here is

\begin{theorem}
\lab{thdiamjie}
(a) (upper bound)  Let $M$ be a compact Riemannian manifold of dimension $n \ge 3$.
Let $g=g(t)$, $t \in [0, T]$ be a family of smooth metrics evolving
under the Ricci flow $\partial_t g = - 2 Ric$, or a static metric. In the later case the time variable is mute. Then there exists
a positive constant $C$, depending only on the Sobolev constants $A$ and $B$,
or the positive Yamabe constant $Y(g)$,  such that
\[
diameter (M, g(t)) \le C \left( 1+ Volume (M, g(t)) +
\int_M R^{(n-1)/2}_{+}(x, t) dg(t) \right).
\] Here $R_+$ is the positive part of the scalar curvature under $g(t)$.

(b) (lower bound.)   Let $(M, g(t))$ be a Ricci flow given in (a).
Then either $diam(M, g(t)) \ge \sqrt{t}$ or
\[
diam(M, g(t))  \ge H(t, g(0)) \,  [V(M, g(0))]^{\frac{1}{n}} \, e^{-
\frac{2}{n}
 \int^{t}_0 \Vert R(\cdot, s) \Vert_{\infty} ds}.
\] where $H$ is a positive function defined by
\[
H= c e^{ \frac{1}{n} [-\alpha - t \beta- t \, \Vert R_-(\cdot, 0)
\Vert_\infty]}
 [ 1+ \frac{2}{n} \Vert R_-(\cdot, 0) \Vert_\infty t
]^{-1/2}.
\] Here $c$ is an absolute constant; $\alpha$ and
  $\beta$ are positive constants depending only on the Sobolev constants of $({\M}, g(0))$ and the
  infimum of Perelman's $F$ entropy for $(M, g(0))$. Also
  $\beta=0$ if $R(\cdot, 0) \ge 0$.
\end{theorem}

{\remark
\label{re}

(1).  It will be shown in the proof that $C=C_0 (A+B+1)^{n/2}$ where
$C_0$ is a constant depending only on $n$. Alternatively $C=C_0 Y(g)^{-n/2}$.

(2). In the case of Ricci flows, if the infimum of the $F$ entropy is positive, then $C$ is
independent of time.  These will be pointed out during
the proof.
Applying the maximum principle on the equation for the scalar curvature
\[
\Delta R - \partial_t R + 2 |Ric|^2=0,
\]it is known that
\[
Volume (M, g(t)) \le Volume (M, g(0)) \, [ \frac{2}{n} \Vert R_-(\cdot, 0) \Vert_\infty  t
+1 ]^{n/2}.
\]Here $R_-= - \min\{0, R\}.$

(3). As mentioned earlier, the parameter $\beta=0$ if the
  initial scalar curvature is nonnegative. In this case we get the
  lower
  bound
\[
diam(M, g(t))\\
 \ge c e^{ - \frac{1}{n} \alpha } e^{- \frac{2}{n}
 \int^{t}_0 \Vert R(\cdot, s) \Vert_\infty ds} [V(M, g(0))]^{1/n}.
\] This dependence on the scalar curvature is sharp in general as can be
seen through the standard Ricci flow on $S^n$. There the scalar
curvature is $R=C/(T-t)$ and the diameter shrinks to $0$ at time
$T$. }

The following are some notations to be frequently used.  $B(x, r, t)$ denotes
the geodesic ball, centered at $x$, with radius $r$, under the metric
$g(t)$; $|B(x, r, t)|$ is the volume of the ball under $g(t)$;  if $\Gamma$ is a curve, then $|\Gamma|$ denotes its length under $g(t)$; $R$ is the scalar curvature. If no confusion arises, we will suppress the time variable $t$.
\medskip
\section{proof of theorem}

\noindent{\it Proof of the theorem, part (a),  the upper bound.}

We will just deal with the case where the Sobolev constants are involved.
The one for Yamabe constant can be treated in exactly the same way.
 Let $x$ be a point in $M$ and $r$ be a positive number less than $diam(M,
g(t))/2$. Following \cite{Ak:1} and \cite{Ca:1}, we take $v=v(y)=r-d(x, y)$ in
(\ref{sob}), or (\ref{sobRF}). Here $d(x, y)$ is the distance between
$x$ and $y$ under the metric $g(t)$. Here and later, if no confusion arises, we suppress the time variable $t$
for brevity of presentation.  Then
\be
\lab{intDv2}
\int_M  4 | \nabla v|^2 dg = 4 |B(x, r)|.
\ee
\be
\lab{intRv2}
\al
\int_M R v^2 dg &\le r^2 \int_{B(x, r)} R_+ dg
=\frac{r^2}{|B(x, r)|} \int_{B(x, r)} R_+ dg \, |B(x, r)|\\
&\le M_2(x, t, R_+, r) |B(x, r)|.
\eal
\ee Here $M_2(x, t, R_+, r)$ is a maximal type function defined by, following
\cite{T:1},
\be
\lab{defM2}
M_2(x, t, R_+, r) = \sup_{0<\rho \le r} \frac{\rho^2}{|B(x, \rho, t)|}
\int_{B(x, \rho, t)} R_+ dg(t).
\ee Also
\be
\lab{intv2}
 \int_{B(x, r)}  v^2 dg \le r^2 |B(x, r)|.
\ee

In the ball $B(x, r/2)$, it is clear that $v \ge r/2$. Therefore, after using
H\"older inequality and (\ref{sob}), we obtain
\[
\al
\frac{r^2}{4} |B(x, r/2)| &\le \int_{B(x, r)}  v^2 dg\\
&\le  |B(x, r)|^{2/n} \, \left(\int_{B(x, r)} v^{2n/(n-2)} dg \right)^{(n-2)/n} \\
& \le |B(x, r)|^{2/n} \left[ A \int_M ( 4 | \nabla v|^2 + R v^2 ) dg(t) + B \int_M v^2 dg(t) \right].
\eal
\] Substituting (\ref{intDv2}),  (\ref{intRv2}) and  (\ref{intv2}) into the
right hand side of this inequality, we find that
\[
\frac{r^2}{4} |B(x, r/2)| \le |B(x, r)|^{(n+2)/n} [ 4 A + 4 A M_2(x,
t, R_+, r) + B r^2].
\]This implies
\[
|B(x, r)| \ge |B(x, r/2)|^{n/(n+2)}  r^{2 n/(n+2)}
[ 16 A (1+ M_2(x, t, R_+, r) )   + 4 B r^2 ]^{-n/(n+2)}.
\]For any number $s \in (0, r]$, it is obvious that the above inequality still
holds when $r$ is replaced by $s$. Since
\[
M_2(x, t, R_+, s)
\le M_2(x, t, R_+, r), \qquad  s^2 \le r^2,
\]we arrive at the following inequality for all $s \in (0, r]$,
\be
\lab{BxsBxs/2}
|B(x, s)| \ge |B(x, s/2)|^{n/(n+2)}  s^{2 n/(n+2)}
[ 16 A (1+ M_2(x, t, R_+, r) )   + 4 B r^2 ]^{-n/(n+2)}.
\ee

Iterating (\ref{BxsBxs/2}) with $s=r, r/2, ..., r/2^m$ for positive integers $m$,
we deduce
\[
|B(x, r)| \ge \left[ r^2 2^{-2} \left( 16 A (1+ M_2(x, t, R_+, r) )   + 4 B r^2\right)^{-1}
\right]^{\Sigma^m_{i=1} (n/(n+2))^i} \, |B(x, r/2^m)|^{(n/(n+2))^m}.
\]Letting $m \to \infty$, this shows, for $r < diam (M, g(t))/2$,
\be
\lab{Qkappa}
\frac{|B(x, r)|}{r^n} \ge
\left[ 64 A (1+ M_2(x, t, R_+, r) )   + 16 B r^2\right]^{-n/2}.
\ee In the Ricci flow case, this can be regarded as a quantified version of Perelman's
$\kappa$ non-
collapsing theorem.

\medskip

{\it Step 2.} Now we start to bound the diameter of $(M, g(t))$, $t \in (0, T]$.
We use $Z=Z(t)$ and $V=V(t)$ to denote the diameter and volume of
$(M, g(t))$ respectively. If no confusion arises, we will ignore $t$.
Without loss of generality, we assume $Z \ge 2$. From now on we assume
$r \le 1$ in (\ref{Qkappa}) so that $B(x, r)$ is always a proper ball in $M$.
Picking any $x \in M$ and writing
\be
\lab{defkap}
\kappa =\kappa(x, r) =
 \frac{|B(x, r)|}{r^n}
\ee in
(\ref{Qkappa}), we see that
\[
64 A (1+ M_2(x, t, R_+, r) )   + 16 B r^2 \ge \kappa^{-2/n}.
\]This implies, as $r \le 1$, that
\be
\lab{M2>}
M_2(x, t, R_+, r)  \ge (64 A)^{-1} (\kappa^{-2/n}-64 A - 16 B).
\ee Set
\be
\lab{defkap0}
\kappa_0 = \min \{(128 A + 16 B)^{-n/2}, \, \frac{\omega_n}{2} \},
\ee where $\omega_n$ is the volume of $n$ dimensional Euclidean ball,
which is here for later use. From (\ref{M2>})  and (\ref{defkap0}), we
know that the following statement is true: if $\kappa = \frac{|B(x, r)|}{r^n}
\le \kappa_0$, then
\be
\lab{M2>2}
M_2(x, t, R_+, r) \ge 2.
\ee
\medskip

{\it Step 3.} Recall that $Z$ is the diameter of $(M, g(t))$. Let $N$ be the greatest integer which is less than or equal to
$Z/4$. Let $a$ and $b$ be two points in $M$ such that $d(a, b)= Z$ . Let $\Gamma$ be a minimum geodesic connecting $a$ and $b$. Let $p$ be the
middle point of $\Gamma$ so that $d(a, p)=d(b, p) = Z/2$.

 Next we claim that if
\be
\lab{Z>}
Z > \frac{V 4^{n+3}}{\kappa_0},
\ee then for at least $N$ many positive integers $i \le
2 N$, it holds
\[
|B(p, i)-B(p, i-1)| \le \kappa_0 4^{-n}.
\] The proof is simple, for, if the claim were not true, then there would be
at least $N$ many $i$ such that
\[
|B(p, i)-B(p, i-1)| > \kappa_0 4^{-n}.
\] which would imply
\[
Z \kappa_0 4^{-n-2} \le N \kappa_0 4^{-n} \le \Sigma^{ 2 N}_{i=1} |B(p, i)-B(p, i-1)| \le V.
\] Thus
\[
Z \le \frac{ V 4^{n+2}}{\kappa_0}.
\]
But this contradicts with (\ref{Z>}), proving the claim. From now
on, we always assume that (\ref{Z>}) holds. This does not reduce any
generality since the theorem is already proven otherwise.

According to the claim, we can pick $N$ integers $i_1, ..., i_N$ in the set $\{1, ..., 2N\}$ such that
\be
\lab{B-B<}
|B(p, i_j)-B(p, i_j-1)| \le \kappa_0 4^{-n}, \qquad j=1, ..., N.
\ee Pick $i \in \{i_1, ..., i_N\}$. Denote by $\Gamma_i$ the segment
$\Gamma \cap (B(p, i)-B(p, i-1))$. Let $p_i$ be the middle point of $\gamma_i$. Then
\[
B(p_i, 1/2) \subset B(p, i) - B(p, i-1).
\]Hence, for any $x \in B(p_i, 1/4)$, we have
\[
B(x, 1/4) \subset B(p_i, 1/2) \subset B(p, i) - B(p, i-1).
\] This and (\ref{B-B<}) infer that
\be \lab{B1/4<} \frac{| B(x, 1/4) |}{(1/4)^n}  \le \kappa_0. \ee On
the other hand
\[
\lim_{\rho \to 0} \frac{|B(x, \rho)|}{\rho^n} = \omega_n \ge 2 \kappa_0,
\]where the last inequality is due to the definition of $\kappa_0$ in
(\ref{defkap0}). From this and (\ref{B1/4<}), we can find a positive number
$s=s(x) \in (0, 1/4]$, which satisfies the following properties.
First
\be
\lab{Bsx=k0}
\frac{| B(x, s(x)) |}{s(x)^n}  = \kappa_0.
\ee Second, if $0<\rho \le s(x)$, then
\be
\lab{Brho>k0}
\frac{|B(x, \rho)|}{\rho^n} \ge \kappa_0.
\ee In another word, $s(x)$ is the smallest radius $\rho$ such that
$\frac{|B(x, \rho)|}{\rho^n} = \kappa_0.$
From (\ref{Bsx=k0}) and
(\ref{M2>2}) with $r=s(x)$, we find that
\[
M_2(x, t, R_+, s(x)) \ge 2.
\] Hence there exists $s_1(x) \in (0, s(x)]$ such that
\be
\lab{s1B}
\frac{s_1(x)^2}{|B(x, s_1(x))|} \int_{B(x, s_1(x))} R_+ dg \ge 1.
\ee According to (\ref{Brho>k0}), we also have
\be
\lab{Bxs1>}
\frac{|B(x, s_1(x))|}{s_1(x)^n} \ge \kappa_0.
\ee

The rest of the proof is similar to that in \cite{T:1}. It is here for completeness.  Denote by $\sigma$ the (disjointed) curve
\[
\cup_{i=i_1, ..., i_N} \left( \Gamma_i \cap B(p_i, 1/4) \right).
\]Since $N \ge Z/8$, we see that
\[
|\sigma| \ge Z/16.
\]Also the family of balls $\{ B(x, s_1(x)) \, | \, x \in \sigma \}$
forms an open cover of $\sigma$. By standard argument, see Lemma 5.2
in \cite{T:1}, e.g., we can find a sequence of points $\{x_l \, | \,
l=1, 2, ...\} \subset \sigma$ such that each of the balls $B(x_l,
s_1(x_l))$ are disjoint from each other and that the balls $\{B(x_l,
s_1(x_l)) \, | \, l=1, 2, ... \}$ cover at least $1/3$ of $\sigma$.
Consequently
\be
\lab{Z16} Z \le 16 | \sigma | \le 96 \Sigma_l |
s_1(x_l)|. \ee Using (\ref{s1B}) and H\"older inequality, we have
\[
\frac{|B(x, s_1(x_l))|}{s_1^2(x_l)} \le \int_{B(x, s_1(x_l))
} R_+ dg
\le \left( \int_{B(x, s_1(x_l))} R^{(n-1)/2}_+ dg \right)^{2/(n-1)} \,
|B(x,  s_1(x_l))|^{(n-3)/(n-1)}.
\] This and (\ref{Bxs1>}) together imply that
\[
\kappa_0 s_1(x_l) \le \frac{|B(x, s_1(x_l))|}{s_1^{n-1}(x_l)}
\le \int_{B(x, s_1(x_l))} R^{(n-1)/2}_+ dg.
\]Using this and (\ref{Z16}), we have proven that
\[
Z \le 96 \kappa_0^{-1} \int_M R^{(n-1)/2}_+ dg.
\] This proves the theorem part (a) since we have assumed $Z \ge 2$ and $Z
\ge \frac{V 4^{n+2}}{\kappa_0}$ to reach the above bound.
Observe that $\kappa_0$, defined on (\ref{defkap0}),  only depends on the Sobolev constants of
$(M, g(t))$ and dimension $n$, which in turn depend only on $t$ and the Sobolev
constants of $(M, g(0))$. If the infimum of the $F$ entropy is positive then
$\kappa_0$ is independent of time. See Theorem 6.2.1 in \cite{Z11:1}.

Finally, the constant in the theorem, part (a),  is as claimed in Remark \ref{re} due to the size of $\kappa_0$.
\medskip

{\it Proof of Part (b).}

Now we prove the lower bound for the diameter.  Let $G
=G(z, l; x, t)$ be the fundamental solution of the conjugate heat equation, with $l<t$. By (1.3) in \cite{Z12:1}, we have the following lower
bound for $G$.
 \be
 \lab{Gxiajie0}
 G(z, l; x, t) \ge \frac{c_1 J(t)}{ (t-l)^{n/2}}
e^{-2 c_2 \frac{d(z, x, t)^2}{t-l}}  e^{- \frac{1}{ \sqrt{t-l}}
\int^{t}_{l} \sqrt{t-s}
 R(x, s) ds}.
 \ee Here
 \be
 \lab{J=}
 J=J(t)=\exp[-\alpha - t \beta- t \, \sup R_-(\cdot, 0)],
 \ee and $\alpha$ and
  $\beta$ are positive constants depending only on the Sobolev constants of $({\M}, g(0))$
  and the
  infimum of Perelman's $F$ entropy for $(M, g(0))$. Moreover $\beta =0$ if $R(\cdot, 0) \ge 0$. Therefore
\be
 \lab{Gxiajie}
G(z, l; x, t) \ge  \frac{c_1 J(t)}{ (t-l)^{n/2}} e^{-2 c_2
\frac{d(z, x, t)^2}{t-l}}  e^{- \int^{t}_{l}
 \Vert R(\cdot, s) \Vert_\infty ds}.
 \ee

Fix a time $t_0>0$ and a point $x_0 \in M$.  Write $r =
\frac{diam(M, g(t_0))}{2}$.  If $r \ge \sqrt{t_0}$, then we do not need
to do anything. So we assume $r <\sqrt{t_0}$.
 In (\ref{Gxiajie}), we take $z=x_0$, $t=t_0$ and $l=t_0-r^2$.
Thus, for $x$ such that $d(x_0, x, t_0) \le r$, we obtain
 \be
\lab{G>r-n}
 G(x_0, t_0-r^2; x, t_0)  \ge  \frac{c_1 J(t_0)}{r^n}
e^{-2 c_2 - \int^{t_0}_0
 \Vert R(\cdot, s) \Vert_\infty ds}
\ee By simple differentiation in time and applying the maximum principle on the scalar
curvature, it is easy to see that
\[
\Vert R_-(\cdot, t) \Vert_\infty \le \frac{1}{(1/\Vert R_-(\cdot,
0) \Vert_\infty) + (2t/n)} \] and
\[
\frac{d}{d t} \int_{\M} G(z, l; x, t) dg(x, t) \le \Vert
R_-(\cdot, t) \Vert_\infty \int_{\M} G(z, l; x, t) dg(x, t).
\]Therefore
 \be
 \lab{intG}
 \int_{\M} G(z, l; x, t) dg(x, t) \le \left[ 1+ \frac{2}{n} \Vert R_-(\cdot,
0) \Vert_\infty (t-l) \right]^{n/2}.
 \ee  Substituting (\ref{G>r-n}) to (\ref{intG}), we deduce
\[
\al \left[ 1+ \frac{2}{n} \Vert R_-(\cdot,
0) \Vert_\infty r^2 \right]^{n/2}& \ge \int_{\M} G(x_0, t_0-r^2; x, t_0) dg(x, t_0)\\
&\ge \int_{d(x_0, x, t_0) \le r} G(x_0, t_0-r^2; x, t_0) dg(x, t_0) \\
& \ge \frac{c_1 J(t_0)}{r^n} e^{-2 c_2 -
  \int^{t_0}_0
 \Vert R(\cdot, s) \Vert_\infty ds} \int_{d(x_0, x,
t_0) \le r}  dg(x, t_0).
 \eal
\] Since $r =
\frac{diam(M, g(t_0))}{2}$ by choice, we know that
$
\int_{d(x_0, x,
t_0) \le r}  dg(x, t_0) = V(M, g(t_0))$, the volume of $M$.
Notice that
\[
\frac{d}{dt} V(M, g(t)) =- \int_M R(x, t) dg(t) \ge -  \Vert
R(\cdot, t) \Vert_{\infty} V(M, g(t)).
\]Hence
\[
 V(M, g(t_0)) \ge e^{- \int^{t_0}_0 \Vert R(\cdot, t) \Vert_{\infty} dt} V(M, g(0)).
\] These imply that
\[
\left[ 1+ \frac{2}{n} \Vert R_-(\cdot, 0) \Vert_\infty r^2
\right]^{n/2} \, r^n
 \ge c_1 J(t_0) e^{-2c_2- 2
 \int^{t_0}_0 \Vert R(\cdot, t) \Vert_{\infty} dt}  \, V(M, g(0)).
\]Since $r=diam(M, g(t_0))/2$ and $ 2 r \ge \sqrt{t_0}$ by assumption,
we see that
\[
\al &diam(M, g(t_0))\\
& \ge c_3 e^{ \frac{1}{n} [-\alpha - t_0 \beta- t_0 \, \Vert
R_-(\cdot, 0) \Vert_\infty]} e^{- \frac{2}{n}
 \int^{t_0}_0 \Vert R(\cdot, t) \Vert_{\infty} dt}
 [ 1+ \frac{2}{n} \Vert R_-(\cdot, 0) \Vert_\infty t_0
]^{-1/2} [V(M, g(0))]^{1/n}. \eal
\]
Here we just used (\ref{J=}) so that $\alpha$ and
  $\beta$ are positive constants depending only on the Sobolev constants of $({\M}, g(0))$ and the
  infimum of Perelman's $F$ entropy for $(M, g(0))$. Also $c_3$ is
  an absolute constant. As mentioned earlier, $\beta=0$ if the
  initial scalar curvature is nonnegative. In this case we get the
  bound
\[
diam(M, g(t_0))\\
 \ge c_3 e^{ - \frac{1}{n} \alpha } e^{- \frac{2}{n}
 \int^{t_0}_0 \Vert R(\cdot, t) \Vert_\infty dt} [V(M, g(0))]^{1/n}.
\] The proof is complete. \qed

\medskip

{\noindent \it Acknowledgments. We wish to thank Professors Mingliang Cai, Bennett Chow and Peter Topping for useful
suggestions about the presentation of the note.}

\bigskip

\noindent e-mail:  qizhang@math.ucr.edu


\begin{thebibliography}{00}

\bibitem [Ak]{Ak:1} Akutagawa, Kazuo, {\it Yamabe metrics of positive scalar
curvature and conformally flat manifolds}. Differential Geom. Appl.
4 (1994), no. 3, 239--258.

\bibitem[Ca]{Ca:1} Carron, Gilles, {\it
In\'egalit\'es isop\'erim\'etriques de Faber-Krahn et
cons\'equences.} (French) Actes de la Table Ronde de G\'eom\'etrie
Diff\'erentielle (Luminy, 1992), 205--232, S\'emin. Congr., 1, Soc.
Math. France, Paris, 1996

 \bibitem [CG]{CG:1} Mingliang Cai and Gregory J Galloway,
{On the topology and area of higher-dimensional black holes}, Class. Quantum Grav. 18 (2001), 2707-2718.


\bibitem [C++]{C++:1} Chow, Bennett; Chu, Sun-Chin; Glickenstein, David; Guenther, Christine; Isenberg, James; Ivey, Tom; Knopf, Dan; Lu, Peng; Luo, Feng; Ni, Lei,
{\it  The Ricci flow: techniques and applications. Part III}. Geometric-analytic aspects. Mathematical Surveys and Monographs, 163. American Mathematical Society, Providence, RI, 2010.

\bibitem[IK]{IK:1} Ilmanen, Tom; Knopf, Dan, {\it A lower bound for the diameter of
 solutions to the Ricci flow with nonzero $H^1(M^n, R)$.} Math. Res. Lett. 10 (2003), no. 2-3, 161-168.

\bibitem[ST]{ST:1} Sesum, Natasa; Tian, Gang, {\it Bounding scalar curvature and diameter along the K\"ahler Ricci flow (after Perelman).} J. Inst. Math. Jussieu 7 (2008), no. 3, 575-587.


\bibitem[T]{T:1} Topping, Peter, { Diameter control under Ricci flow},
Comm. Ana. Geo., Vol. 13, (2005) 1039-1055.

\bibitem[Z11]{Z11:1}
Qi S. Zhang.
\newblock {\em Sobolev inequalities, heat kernels under Ricci flow and the
  Poincar\'e conjecture}.
\newblock CRC Press, Boca Raton, FL, 2011.

\bibitem[Z12]{Z12:1}
Qi S. Zhang, {\it
Bounds on volume growth of geodesic balls under Ricci flow,}
Math. Res. Letters, 19 (2012), no. 1, 245-253;
    arXiv:1107.4262

\end{thebibliography}
\enddocument